\documentclass{article}
\usepackage{amsthm}
\usepackage[utf8]{inputenc}
\usepackage[T1]{fontenc}
\usepackage[centertags]{amsmath}
\usepackage{amsfonts}
\usepackage{graphicx}
\usepackage{amsmath}
\usepackage{amssymb}
\usepackage{latexsym}
\usepackage[mathscr]{euscript}
\usepackage{tikz,color,makeidx} 
\usetikzlibrary{calc,matrix,arrows,chains,positioning,scopes,decorations.pathmorphing,backgrounds,fit}

\numberwithin{equation}{section}
\newtheorem{teo}{Theorem}[section]
\newtheorem{cor}{Corollary}[section]
\newtheorem{pro}{Proposition}[section]
\newtheorem{lem}{Lemma}[section]
\newtheorem{es}{\textbf{Example}}[section]
\newtheorem{defi}{\mbox{\textbf{Definition}}}[section]
\newtheorem{rem}{\mbox{\textbf{Remark}}}[section]
\newtheorem{notation}{\textit{Notation}}[section]

\newcommand{\bdfn}{\begin{defi} \begin{rm}}
\newcommand{\edfn}{\end{rm} \end{defi}}
\newcommand{\bthm}{\begin{teo}}
\newcommand{\ethm}{\end{teo}}
\newcommand{\bprop}{\begin{pro}}
\newcommand{\eprop}{\end{pro}}
\newcommand{\bcor}{\begin{cor}}
\newcommand{\ecor}{\end{cor}}
\newcommand{\blem}{\begin{lem}}
\newcommand{\elem}{\end{lem}}
\newcommand{\bfact}{\begin{rem} \begin{rm}}
\newcommand{\efact}{\end{rm} \end{rem}}
\newcommand{\bex}{\begin{es} \begin{rm}}
\newcommand{\eex}{ \end{rm} \end{es}}
\newcommand{\ten}{\otimes}
\newcommand{\vsp}[1]{\vspace*{#1cm}}
\newcommand{\bnot}{\begin{notation} \begin{rm}}
\newcommand{\enot}{\end{rm} \end{notation}}
\newcommand{\quot}[2]{{\raisebox{.2em}{$#1$}\left/\raisebox{-.2em}{$#2$}\right.}}

\tikzset{node distance=2.3cm, auto}

\title{Scalar extensions for algebraic structures of \L ukasiewicz logic}

\author{Serafina Lapenta\\ 
 {\small Department of Mathematics, Computer Science and Economics,}\\ 
{\small University of Basilicata,
 Viale dell'Ateneo Lucano, 10 C.P. 85100, Potenza, Italy}\\
 {\small serafina.lapenta@unibas.it}
\and 
Ioana Leu\c stean \\
{\small Department of Computer Science,} \\
{\small Faculty of Mathematics and Computer Science, University of Bucharest,}\\
{\small Academiei nr.14, sector 1, C.P. 010014,  Bucharest, Romania}\\   
{\small ioana@fmi.unibuc.ro}
}

\date{}

\begin{document}
\maketitle
\begin{abstract}
In this paper we study the tensor product  for  MV-algebras, the algebraic structures  of \L ukasiewicz $\infty$-valued logic. Our main results are: the proof that the tensor product is preserved by the categorical equivalence between the MV-algebras and abelian  lattice-order groups with strong unit and the proof of the scalar extension property for semisimple MV-algebras. We explore consequences of this results for various classes of MV-algebras and lattice-ordered groups enriched with a product operation. 
\end{abstract}

\section*{Introduction}

\noindent MV-algebras are the algebraic counterpart of \L ukasiewicz $\infty$-valued logic \cite{Cha1}. 
 The  variety of MV-algebras is generated by 
the standard  model $([0,1], \oplus, ^*, 0)$, where $[0,1]$ is the unit real interval, $x\oplus y=min\{1,x+y\}$ and $x^*=1-x$ for any $x,y\in [0,1]$.  Since the interval $[0,1]$ is closed to the real product, a natural problem  was to  analyze the systems obtained by enriching \L ukasiewicz logic with a product operation. This investigation lead to fruitful investigations both in algebra and logic. 

Several extension of the notion have been defined by endowing an MV-algebra with products: an internal binary product leads to the notion of PMV-algebra \cite{DiND}; a scalar product leads to the notion of MV-module \cite{LeuMod} and Riesz MV-algebra \cite{LeuRMV};  a combination of both leads to the notion of \textit{f}MV-algebra \cite{LL1}. For all these structures, corresponding logical systems are developed. 

 Within several important results, one main achievement in the theory is the categorical equivalence with abelian lattice-ordered groups with strong unit \cite{Mun1}.  The categorical equivalence for MV-algebras extends naturally extended to any MV-algebra with product, and allows us to connect PMV-algebras, MV-modules, Riesz MV-algebras  and \textit{f}MV-algebras with $\ell$-rings with strong unit, $\ell$-modules with strong unit, Riesz Spaces with strong unit and \textit{f}-algebras with strong unit respectively. 

The tensor product and its uses are well-known, therefore its definition was given  by Martinez  for lattice-ordered groups \cite{Mart} and Mundici for MV-algebras  \cite{Mun}.  An important subclass of MV-algebras are the semisimple ones which correspond, through the above mentioned categorical equivalence, to archimedean lattice-ordered structures. Since the tensor product  does not preserve the semisimplicity, the semisimple tensor product was defined in \cite{Mun} for MV-algebras and \cite{Fre,BVanR} for lattice-ordered groups.

The scalar extension property (SEP) is one of the basic properties arising from a tensor product; while it is straightforward in the non-ordered case,  it presents some difficulties in the framework of  lattice-ordered structures. Note that MV-algebras have a natural order wich makes them lattice-ordered structures as well.  In Section \ref{sec:SEP} we give an account of the results known so far, we prove SEP  for the semisimple tensor product and we analyze some of its consequences. 
Since in  Section \ref{sec:gamma} we prove that the tensor product is preserved by the categorical equivalence  between MV-algebras and lattice-ordered groups,  the results are stated  both in the theory of MV-algebras  and  in the theory of lattice-ordered groups. 

The scalar extension property led us to categorical adjunctions between semisimple  MV-algebras and semisimple Riesz MV-algebras  in Section \ref{sec:SEP} and between semisimple PMV-algebras and semisimple \textit{f}MV-algebras in Section \ref{sec:PMV}. In Section \ref{sec:final} we sum up our results and provide an insight on their significance for \L ukasiewicz  logic.

\section{Preliminaries}

\subsection{MV-algebras and $\ell u$-groups}
\noindent An {\em MV-algebra} is an algebraic structure $(A, \oplus, \ ^*, 0)$, where $(A, \oplus ,0)$ is a commutative monoid, $\ ^*$ is an involution and the identity $(x^*\oplus y)^*\oplus y=( y^*\oplus x)^*\oplus x$ is satisfied for any $x,y,z\in A$. We further define $x\odot y = (x^* \oplus y^*)^*$ and $1=0^*$. An order can be defined on $A$ by setting $x\le y$ if and only if $x^*\oplus y=0$; if we set $x\vee y=(x^*\oplus y)^*\oplus y$  and $x\wedge y=(x^*\vee y^*)^*$ then $(A,\vee, \wedge, 0,1)$ is a bounded distributive lattice. 
We refer to \cite{CDM, MunBook} for all the unexplained notions concerning MV-algebras. 

 If $A$ is an MV-algebra we define a partial operation $+$ as follows: for any $x,y\in A$, $x+y$ is defined if and only if $x\leq y^*$ and, in this case, $x+y=x\oplus y$. Note that this operation is cancellative.  Assume $A$, $B$ and $C$ are MV-algebras.  A  function $\omega:A\to B$ is  {\em linear}  if $a\leq b^*$ implies $\omega(a)\leq\omega(b)^*$ and  $\omega(a+b)=\omega(a)+\omega(b)$. Bilinear functions $\beta:A \times B\to C$ are defined as usual.

Semisimple MV-algebras will play an important role in our development. If $A$ is an MV-algebra and  $Rad(A)$ is the intersection of its maximal ideals, then $A$ is semisimple if and only if $A$ is isomorphic to a separating MV-algebra of $[0,1]$-valued continuous functions defined over some compact Hausdorff space \cite{CDM}.

An {\em $\ell u$-group}  is a pair $(G,u)$, where $G$  is an abelian lattice-ordered  group \cite{Birk, BKW} and $u$ is a strong unit. 
If $(G,u)$ is an $\ell u$-group, then $[0,u]_G=([0,u],\oplus, ^*, 0)$ is an MV-algebra, where  $[0,u]=\{ x\in G \mid 0\le x\le u\}$ and  $ x\oplus y=u\wedge (x+y)$, $x^*=u-x$ for any $x\in [0,u]$. 

If $\mathbf{MV}$ is the category of MV-algebras and $\mathbf{auG}$ is the category of $\ell u$-groups equipped with morphisms that preserve the strong unit, then one defines a functor $\Gamma: \mathbf{auG}\rightarrow \mathbf{MV}$ by
$\Gamma (G,u)=[0,u]_G $ and $\Gamma(h)= h|_{[0,u_1]_{G_1}}$,
where  $(G,u)$ is an $\ell u$-group   and   $h: G_1\to G_2$ is a morphism in $\mathbf{auG}$ between $(G_1,u_1)$ and $(G_2, u_2)$.  The functor $\Gamma$ establishes a categorical equivalence between  $\mathbf{auG}$ and $\mathbf{MV}$ \cite{Mun1}. Moreover, through $\Gamma$, semisimple MV-algebras correspond to archimedean $\ell u$-groups. 

In the following $\Gamma({\mathbb R}, 1)$ will be simply denoted $[0,1]$ (depending on context, the MV-algebra structure will be tacitly assumed).  By Chang's completeness theorem \cite{Cha2}, the variety of MV-algebras is generated by $[0,1]$.

\subsection{MV-algebras endowed with a product operation}

Product MV-algebras (PMV-algebras for short) have been defined in \cite{DiND} in the general case and in \cite{MonPMV}  an equivalent axiomatization was provide for the unital and commutative structures.    A {\em unital PMV-algebra} is a structure 
$(P,\oplus,\cdot, ^*, 0)$ such that $(P,\oplus, ^*, 0)$ is an MV-algebra and $\cdot:P\times P\to P$ is a bilinear function  such that 
$a \cdot (b\cdot c)= (a \cdot b)\cdot c$ and $a\cdot 1=1\cdot a=a$  for any $a,b,c\in P$.

A further extension of the notion of MV-algebra has been introduced in \cite{LeuMod}. If $P$ is a PMV-algebra, then 
an  {\em MV-module} over $P$ ({\em $P$-MV-module}) is  a structure $(M,\oplus, ^*,\{\alpha|\alpha\in P\}, 0)$ such that $(M,\oplus, ^*, 0)$ is an 
MV-algebra, $\{\alpha|\alpha\in P\}$ is a family of unary operations  such that the function $(\alpha,x)\mapsto \alpha x$  is bilinear,  $(\alpha \cdot \beta)x= \alpha (\beta x)$ and $1x=x$  for any $\alpha,\beta\in P$ and any $x\in M$.

Note that  \cite[Section 6.4]{LeuMod}  provides an equational characterization for these structures. Most important for our development is the case $P=[0,1]$. The MV-modules over $[0,1]$ are called {\em Riesz MV-algebras}   are studied in  \cite{LeuRMV}. 

Finally, {\em unital \textit{f}MV-algebras} have been introduced in \cite{LL1} and they are algebraic structures $(A, \oplus, \ ^*, \cdot, \{ \alpha \}_{\alpha \in [0,1]} , 0)$ such that $(A, \oplus, \ ^*, \cdot ,0)$ is a unital PMV-algebra, $(A, \oplus, \ ^*, \{ \alpha \}_{\alpha \in [0,1]} , 0)$ is a Riesz MV-algebra and the condition $\alpha (x\cdot y)=(\alpha x)\cdot y=x\cdot (\alpha y)$ is satisfied for any $\alpha \in [0,1]$ and $x,y\in A$.

We defined a hierarchy of algebraic structures, all of them having an MV-algebra reduct.  Hence there are   forgetful functors from 
the categories $\mathbf{PMV}$ of PMV-algebras, $\mathbf{RMV}$ of Riesz MV-algebras and $\mathbf{fMV}$ of \textit{f}MV-algebras  to $\mathbf{MV}$.   For each of this structures one can prove a categorical equivalence with an appropriate class of unital lattice-ordered structures having a lattice-ordered group reduct with a strong unit \cite{DiND, LeuRMV, LL1}.  If $\mathbf{uR}$  is the category of unital $f$- rings with strong unit ($fu$-rings), $\mathbf{uRS}$ is the category of Riesz spaces with strong unit  and $\mathbf{fuAlg}$ is the category of unital \textit{f}-algebras with strong unit ($fu$-algebras), then  the  categorical equivalence are presented in  the following diagram, in which all horizontal arrows are suitable forgetful functors.

\begin{center}
\begin{tikzpicture}
  \node (A) {$\mathbf{uR}$};
  \node (B) [below of=A] {$\mathbf{PMV}$};
  \node (C) [right of=A] {$\mathbf{auG}$};
  \node (D) [below of=C] {$\mathbf{MV}$};
  \node (E) [right of=C] {$\mathbf{uRS}$};
  \node (F) [below of=E] {$\mathbf{RMV}$};
\node(H)[left of=A] {$\mathbf{fuAlg}$};
\node(G)[below of=H] {$\mathbf{fMV}$};
\node(J)[right of=E] {$\mathbf{fuAlg}$};
\node(I)[below of=J] {$\mathbf{fMV}$};
 \draw[->] (A) to node [swap] {$\Gamma_{(\cdot)}$} (B);
\draw[->] (H) to node [swap] {$\Gamma_{f}$} (G);
   \draw[->] (J) to node[swap] {$\mathcal{U}_{(\cdot \ell)}$} (E);
 \draw[->] (A) to node {$\mathcal{U}_{(\cdot \ell)}$} (C);
  \draw[->] (C) to node [swap] {$\Gamma$} (D);
  \draw[->] (B) to node [swap] {$\mathcal{U}_{(\cdot)}$} (D);
\draw[->] (I) to node {$\mathcal{U}_{(\cdot)}$} (F);  
\draw[->] (E) to node [swap] {$\mathcal{U}_{( \ell \mathbb{R})}$} (C);
\draw[->] (H) to node {$\mathcal{U}_{( \ell \mathbb{R})}$} (A);  
\draw[->] (E) to node {$\Gamma_{\mathbb{R}}$} (F);
\draw[->] (G) to node[swap] {$\mathcal{U}_{\mathbb{R}}$} (B);
    \draw[->] (J) to node {$\Gamma_{f}$} (I);
  \draw[->] (F) to node {$\mathcal{U}_{\mathbb{R}}$} (D);
 \node [below=0.5cm, align=flush center] at (D){ Figure 1.};
\end{tikzpicture}
\end{center}

Note that for the objects in $\mathbf{uR}$  and $\mathbf{fuAlg}$ have a  unital ring reduct and we ask, in addition, that the ring-unity coincides with the strong unit.  For the general theory of $\ell$-rings and \textit{f}-rings, Riesz spaces, $\ell$-algebras and \textit{f}-algebras we refer to \cite{BKW, Birk,BP,BVanR,RSZan}. 

We note that a PMV-algebra (Riesz MV-algebra, \textit{f}MV-algebra) is semisimple if its MV-algebra reduct is  a semisimple MV-algebra and that  semisimple PMV-algebras (semisimple Riesz MV-algebras, semisimple \textit{f}MV-algebras) correspond to archimedean $\ell u$-rings (archimedean unital Riesz spaces, archimedean \textit{f}$u$-algebras).  Also, any unital and semisimple PMV-algebra (\textit{f}MV-algebra) is commutative by the general theory of \textit{f}-rings (\textit{f}-algebras) \cite{BP, RSZan}.

In a similar manner, the MV-modules are categorically equivalent to appropriate classes of $\ell$-modules  with strong unit ($\ell u$-rings). We refer to \cite{Stein} for the general theory of $\ell$-modules.  If $(R,u)$ is an $\ell u$-ring and $P\simeq\Gamma_{(\cdot)}(R,u)$ then the category  of MV-modules over $P$ is equivalent to the category of $\ell u$-modules over $R$ \cite{LeuMod}.

\subsection{The  tensor product $\ten_{MV}$ and the semisimple tensor product $\ten$}

 A bimorphism is a bilinear function that is $\vee$-preserving and $\wedge$-preserving in each component. Bimorphisms were defined in  \cite{Mun}, where  the additional requirement $\beta(1, 1)=1$ was imposed. In the present approach we  eliminate this restriction.

The \textit{interval algebra} of $A$ is the MV-algebra $[0,a]=\{ x\in A \mid 0\le x\le a \},$ endowed with the following operations $x\oplus _a y= (x\oplus y)\wedge a$, $x^{*a}=x^* \odot a$ for any $x,y \in [0,a]$ \cite{Mun}. In the follow we will use the notation $[0,a]\le_i A$ in order to say that $[0,a]$ is an interval algebra of $A$.

The MV-algebraic tensor product was defined in \cite{Mun} as a universal bimorphism. We shall use in the sequel a slightly modified universal property proved in  \cite{LeuTens}.   For two MV-algebras $A$ and $B$, let  $A\otimes_{MV} B$ be the tensor product and $\beta_{A,B}:A\times B\to A\otimes_{MV} B$ the universal bimorphism. Then  the  following universal property  holds:

{\em for any MV-algebra $C$ and for any bimorphism $\beta: A\times B\rightarrow C$, there is a unique homomorphism of MV-algebras $\omega :A\otimes_{MV}B\rightarrow [0, \beta (1,1)]\le_i C$ such that $\omega \circ \beta_{A,B}=\beta$.}

\noindent For $a\in A$ and $B\in B$ we denote $a\ten_{MV}b=\beta_{A,B}(a,b)$. Note that $A\ten_{MV}B $ is generated by 
$\beta_{A,B}(A\times B)$. 

In \cite{Mun} the author proves that there exists a semisimple MV-algebra $A$ such that $A \ten_{MV}A$ is not semisimple. Therefore he defines the semisimple tensor product of $A$ and $B$, semisimple MV-algebras, by 
$$ A\ten B = \quot{A\ten_{MV}B }{ Rad(A\ten_{MV}B)}.$$ 
\noindent For the semisimple tensor product the universal property holds  with respect to semisimple  MV-algebras. 

The following representation theorem is crucial for our development. 

\bthm \cite{Mun}  \label{rem:TPMun}
Let $A$ and $B$ be semisimple MV-algebras, and let $X,Y$ be the set such that $A\subseteq C(X)$ and $B\subseteq C(Y)$. Let $\gamma: A\times B \rightarrow C(X\times Y)$ be the map defined by $\gamma (a,b)(x,y)=a(x)b(y)$. Then $\gamma$ is a bimorphism and $A\ten B$ is isomorphic to the MV-subalgebra of $C(X\times Y)$ generated by $\gamma(a,b)$, with $a\in A$ and $b\in B$.
\ethm

Few types of tensor products are defined in the literature of partially-ordered and lattice-ordered groups \cite{Mart, Stein}. We recall the one that will be used in the sequel. 

Let  $G,\ H$ and $L$ abelian lattice-ordered groups. An $\ell$\textit{-bilinear} function is a map $\gamma: G\times H\rightarrow L$ such that $\gamma (x, \cdot)$ and $\gamma (\cdot , y)$ are homomorphisms of $\ell$-groups when $x$ and $y$ are positive elements in $G$ and $H$, respectively. Hence, the tensor product   $G\otimes_{\ell} H$ and its universal $\ell$-bilinear function
$\gamma_{G,H}:G\times H\to G\otimes_{\ell} H$ satisfy a universal property with respect to abelian lattice-ordered  groups and $\ell$-bilinear functions.  For $x\in G$ and $y\in H$ we denote $x\ten_{\ell}y=\gamma_{G,H}(x,y)$. 

In \cite{BVanR, Fre} the authors provide a construction for the tensor product of archimedean $\ell$-groups, denoted by $ \ten_{a}$ and  they prove the universal property with respect to archimedean structures.

\section{The $\Gamma$-functor and the tensor product $\ten_{MV}$} \label{sec:gamma}

In this section we prove that the functor $\Gamma: \mathbf{auG}\rightarrow \mathbf{MV}$ commutes with  the tensor product, both in the general and in the archimedean case.  To do this, we prove an extension result for bimorphisms. 

As a preliminary step, we give the detailed proof of the fact that any factor of the MV-algebraic  tensor product is embedded  in the tensor product. If $A$ and $B$ are MV-algebras, then we define
\vsp{-0.2}
\begin{center}
 $\iota_A: A \rightarrow A\ten_{MV}B$ and $\iota_{B}: B \rightarrow A\ten_{MV}B$

 $\iota_A(a)=a\ten_{MV}1_B$ and $\iota_B(b)=1_A\ten_{MV}b$ for any $a\in A$ and $b\in B$.
\end{center}
\vsp{-0.2}
The functions $\iota_A$ and $\iota_B$ are  embeddings of MV-algebras \cite{MoFl} (private communication).   We sketch the proof, for the sake of completeness.

\bprop \label{teo:embmvTP}
The maps $\iota_A: A \rightarrow A\ten_{MV}B$ and $\iota_{B}: B \rightarrow A\ten_{MV}B$ defined as $\iota_A(a)=a\ten_{MV}1_B$ and $\iota_B(b)=1_A\ten_{MV}b$ for any $a\in A$ and $b\in B$, are embeddings.
\eprop
\begin{proof}
By \cite[Theorem 2.20]{MunBook} there exists a MV-algebra $D$ such that both $A$ and $B$ embeds in it. By \cite[Theorem 9.5.1]{CDM} there exists a set $X$ and a MV-algebra embedding $f:D \hookrightarrow (^*[0,1])^X$, therefore there exist two embeddings $A\stackrel{f_A}{\hookrightarrow} (^*[0,1])^X$ and $B\stackrel{f_B}{\hookrightarrow} (^*[0,1])^X$. We remark that $(^*[0,1])^X$ is a unital and commutative PMV-algebra, therefore we define the following bimorphism

$\beta : A\times B \rightarrow (^*[0,1])^X$, $\beta_1(a,b)=f_A(a)\cdot f_B(b)$ for any $a\in A$ and any $b\in B$.\\
By the universal property in \cite{Mun}, there exists $\omega: A\ten_{MV}B \rightarrow (^*[0,1])^X$ such that $\omega(a\ten_{MV}b)=f_A(a) \cdot f_B(b)$.\\
Assume that $\iota_A(a_1)=\iota_A(a_2)$, that is $a_1\ten_{MV}1_B=a_2\ten_{MV}1_B$, then $f_A(a_1)=\omega(a_1\ten_{MV}1_B)=\omega(a_2\ten_{MV}1_B)=f_A(a_2)$. Since $f_A$ is an embedding, the conclusion follows. Analogously, $\iota_B$ is an embedding.
\end{proof}

Assume $ (G, u_G)$ $(H,u_H)$ and $(L, u_L)$ are $\ell u$-groups and $\gamma:G\times H\to L$  is an $\ell$- bilinear function. 
We say that $\gamma$ is {\em $\ell u$-bilinear} if $\gamma(u_G, u_H)\leq u_L$. 

In the sequel we prove that a bimorphism uniquely extends to an $\ell u$-bilinear function.  

\bprop\label{lem:01ext}
If  $A=\Gamma (G, u_G),\ B=\Gamma (H,u_H)$ and $C=\Gamma (L, u_L)$  and  $\beta: A\times B \rightarrow C$ is a bimorphism, then there exists an  unique $\ell u$-bilinear function  $\overline{\beta}: G\times H \rightarrow L$ that extends $\beta $. 
\eprop
\begin{proof}
Let $a$ be a fixed element in $A$ and denote $\beta_a=\beta(a,\cdot)$. Hence, by \cite[Proposition 2.3]{Mun},  we have $\beta_a: B\rightarrow [0, \beta (a, 1)]\le_i C$ is a homomorphism of MV-algebras.  
By \cite[Proposition 2.9]{LeuTens}, there exists a unique homomorphism of $\ell$-groups $\overline{\beta_a}$ such that $\overline{\beta _a} :H \rightarrow L$, and $\overline{\beta _a} |_{B}=\beta_a. $

\textit{Step 1}. We prove that the map $ \gamma_h :A \rightarrow L$, defined by $\gamma_h(a)=\overline{\beta_a}(h)$ for any $h\in H$, is linear.\\
Let the sum $a_1+a_2$ be defined in $A$ and let $h\in H$, since $B$ generates the positive cone of $(H, u_H)$, $h=h^+-h^-$, where $h^+=s_1+\ldots +s_n$ and $h^-=t_1+\ldots +t_m$, with $s_i, t_j\in B$ for any $i=1, \ldots , n$ and $j=1, \ldots , m$.
\begin{center}
$ \gamma_h(a_1+a_2)= \overline{\beta_{a_1+a_2}}((s_1+\ldots +s_n)-(t_1+\ldots
 +t_m))=(\overline{\beta_{a_1+a_2}}(s_1)+\ldots + \overline{\beta_{a_1+a_2}
 }(s_n))-(\overline{\beta_{a_1+a_2}}(t_1)+\ldots + \overline{\beta_{a_1+a_2}}
 (t_m))=(\beta(a_1+a_2,s_1)+\ldots +\beta(a_1+a_2,s_n))-(\beta(a_1+a_2,t_1)+\ldots
  \beta(a_1+a_2,t_m))=(\beta(a_1, s_1)+\beta (a_2,s_1)+\ldots +\beta(a_1, s_n)+\beta 
  (a_2,s_n))-(\beta(a_1, t_1)+\beta (a_2,t_1)+\ldots \beta(a_1, t_m)+\beta 
  (a_2,t_m))=\left[ (\overline{\beta_{a_1}}(s_1)+ \ldots +\overline{\beta_{a_1}}
  (s_n))-(\overline{\beta_{a_1}}(t_1)+\ldots +\overline{\beta_{a_1}}(t_m))\right] +
   \left[ (\overline{\beta_{a_2}}(s_1)+ \ldots +\overline{\beta_{a_2}}
   (s_n))-(\overline{\beta_{a_2}}(t_1)+\ldots +\overline{\beta_{a_2}}(t_m))\right] = 
   \overline{\beta_{a_1}}(h)+\overline{\beta_{a_2}}(h)=\gamma_h(a_1)+\gamma_h(a_2).$
\end{center}

\textit{Step 2}. $\gamma_h$ commutes with $\wedge$ and $\vee$.\\
We first remark that an element of $(H,u_H)$ is a good sequence according with Mundici construction of the inverse of the functor $\Gamma$. Moreover, by \cite{LeuGS} and it is possible to take indexes in $\mathbb{Z}$ instead of quotients of sequence, then
\begin{center}
$\gamma_h(a_1\wedge a_2)=\overline{\beta_{a_1\wedge a_2}}(h)=\overline{\beta_{a_1\wedge a_2}}((h_i)_{i\in \mathbb{Z}}) =(\overline{\beta_{a_1\wedge a_2}}(h_i))_{i\in \mathbb{Z}}=(\beta(a_1\wedge a_2, h_i))_{i\in \mathbb{Z}}= (\beta (a_1, h_i)\wedge \beta (a_2, h_i))_{i\in \mathbb{Z}}=(\beta (a_1, h_i))_{i\in \mathbb{Z}}\wedge (\beta (a_2, h_i))_{i \in \mathbb{Z}}=(\overline{\beta_{a_1}}(h_i))_{i\in \mathbb{Z}}\wedge (\overline{\beta_{a_2}}(h_i))_{i\in \mathbb{Z}}=\overline{\beta_{a_1}}(h) \wedge \overline{\beta_{a_2}}(h)= \gamma_h(a_1)\wedge \gamma_h(a_2),$
\end{center}
and similarly for $\vee$.\\
Therefore the map $\gamma_h: [0, u_G]\rightarrow L $ is linear and commutes with $\vee$ and $\wedge$, i.e. $\gamma_h :[0, u_G]\to [0, \gamma_h(u_G)] \le _i [0, u_L]$ is an homomorphism of MV-algebras \cite[Proposition 2.3]{Mun} and consequently, by \cite[Proposition 2.9]{LeuTens} there exists a unique homomorphism of  $\ell$-groups $ \overline{\gamma_h}:G\rightarrow K \le L,$
where $K$ is the $\ell$-group generated by $\gamma_h(A)$  and  $\overline{\gamma_h}|_A=\gamma_h$.\\
We define now $\overline{\beta}:G\times H\rightarrow L$ as $\overline{\beta}(g,h)=\overline{\gamma_h}(g)$. By the construction, $\overline{\beta}(\cdot, h)$ is a homomorphism of $\ell$-groups.

\textit{Step 3}. $\overline{\beta}(g, \cdot)$, with $g$ fixed element in $G$, is linear.\\
Let $h_1, h_2$ be elements in $H$; there exist suitable elements in the unit interval such that $g^+=s_1+\ldots +s_n$ and $g^-=t_1+\ldots +t_m$.
\begin{center}
$ \overline{\beta}(g, h_1+h_2)=\overline{\gamma_{h_1+h_2}}(g)= (\overline{\gamma_{h_1+h_2}}(s_1)+\ldots +\overline{\gamma_{h_1+h_2}}(s_n))-(\overline{\gamma_{h_1+h_2}}(t_1)+\ldots +\overline{\gamma_{h_1+h_2}}(t_m))=(\overline{\beta_{s_1}}(h_1+h_2)+\ldots +\overline{\beta_{s_n}}
(h_1+h_2))-(\overline{\beta_{t_1}}(h_1+h_2) + \ldots +\overline{\beta_{t_m}}(h_1+h_2) 
)= (\overline{\beta_{s_1}}(h_1)+\overline{\beta_{s_1}}(h_2)+\ldots 
+\overline{\beta_{s_n}}(h_1)+\overline{\beta_{s_n}}(h_2))-(\overline{\beta_{t_1}}
(h_1)+\overline{\beta_{t_1}}(h_2) + \ldots +\overline{\beta_{t_m}}
(h_1)+\overline{\beta_{t_m}}(h_2) )=\left[ (\gamma_{h_1}(s_1)+\ldots +\gamma_{h_1}
(s_n))-(\gamma_{h_1}(t_1)+\ldots +\gamma_{h_1}(t_m))\right] + \left[ (\gamma_{h_2}
(s_1)+\ldots +\gamma_{h_2}(s_n))-(\gamma_{h_2}(t_1)+\ldots +\gamma_{h_2}(t_m))\right] =\overline{\gamma_{h_1}}(g)+\overline{\gamma_{h_2}}(g)= \overline{\beta}(g, h_1)+\overline{\beta}(g,h_2).$ 
\end{center}

\textit{Step 4}. $\overline{\beta}(g, \cdot)$ commute with $\vee$ and $\wedge$.\\
We will use again good sequences
\begin{center}
$\overline{\beta}(g, h_1\wedge h_2)=\overline{\gamma_{h_1\wedge h_2}}
(g)=\overline{\gamma_{h_1\wedge h_2}}((g_i)_{i\in 
\mathbb{Z}})=(\overline{\gamma_{h_1\wedge h_2}}(g_i))_{i\in \mathbb{Z}} = 
(\gamma_{h_1\wedge h_2}(g_i))_{i\in \mathbb{Z}} = (\overline{\beta_{g_i}}(h_1\wedge 
h_2))_{i\in \mathbb{Z}}= (\overline{\beta_{g_i}}(h_1)\wedge \overline{\beta_{g_i}}
(h_2))_{i\in \mathbb{Z}}= (\overline{\beta_{g_i}}(h_1))_{i\in \mathbb{Z}}\wedge 
(\overline{\beta_{g_i}}(h_2))_{i\in \mathbb{Z}}=(\overline{\gamma_{h_1}}(g_i))_{i\in 
\mathbf{Z}}\wedge (\overline{\gamma_{h_2}}(g_i))_{i\in \mathbf{Z}}= 
\overline{\gamma_{h_1}}(g) \wedge \overline{\gamma_{h_2}}(g)= \overline{\beta}(g, 
h_1) \wedge \overline{\beta}(g, h_2).
$
\end{center}
The same can be done for $\vee$ and $\overline{\beta}(g, \cdot)$ is a homomorphism of $\ell$-groups. Moreover, $\overline{\beta}(u_G, u_H)= \beta (u_G, u_H)\le u_L$.

In order to prove the uniqueness, we assume that $\tilde{\beta}:G\times H\to L$ is another $\ell u$-bilinear function that extends $\beta$ . If $a\in A$ then $\tilde{\beta}(a,\cdot)$ is an extension of $\beta_a$, so it coincides with 
$\overline{\beta_a}$. It follows that $\overline{\beta}(a,\cdot):H\to L$ and $\tilde{\beta}(a,\cdot):H\to L$ coincide for any $a\in A$. By linearity, they coincide for any $g\in G$.  
\end{proof}

The main result of this section  is Theorem \ref{pro:02ten-iso}, which asserts that the functor $\Gamma: \mathbf{auG}\to \mathbf{MV}$ preserves the tensor product.

Note that, if  $(G, u_G)$ and $(H,u_H)$ are $\ell u$-groups then $u_G \ten _{\ell} u_H$ is strong unit in $G\ten _{\ell} H$  \cite[3.6]{Mart}.  In the sequel we  prove  two preliminary lemmas.

\blem \label{lem:01ten}
Let $(G,u_{G})$, $(H,u_{H})$ and $(L,u_{L})$ be $\ell u$-groups. For any bimorphism
$\gamma:\Gamma(G,u_{G})\times \Gamma(H,u_{H})\rightarrow \Gamma(L,u_{L})$ there is a unique homomorphism of  MV-algebras
$\omega:\Gamma(G\ten_{\ell} H,u_{G}\ten_{\ell}u_{H})\rightarrow [0,\gamma(u_{G},u_{H})]\leq_{i} \Gamma(L,u_{L})$ such that $\omega(x\ten_{\ell}y)=\gamma(x,y)$ for any $x\in \Gamma(G,u_{G})$ and $y\in\Gamma(H,u_{h})$.
\elem
\begin{proof}
 We set 
$A=\Gamma(G,u_{G})$, $B=\Gamma(H,u_{H})$, $C=\Gamma(L, u_{L})$ and we suppose that
  $\gamma:A\times B\rightarrow C$ is a bimorphism.
By Proposition   \ref{lem:01ext}, there is a $\ell u$-bilinear function  $\widetilde{\gamma}:G\times H\rightarrow L$ 
which extends $\gamma$, so there is a unique homomorphism of 
$\ell$-groups $\widetilde{\omega}:G\ten_{\ell}H\rightarrow L$ such that $\widetilde{\omega}\circ\gamma_{G,H}=\widetilde{\gamma}$.\\
It follows that $\widetilde{\omega}(u_{G}\ten_{\ell} u_{H})=\widetilde{\gamma}(u_{G},u_{H})=\gamma(1,1)\leq u_{L},$ therefore by  \cite[Lemma 2.8]{LeuTens} the restriction $\omega:\Gamma(G\ten_{\ell} H,u_{G}\ten_{\ell} u_{H})\rightarrow [0,\gamma(u_{G},u_{H})]\leq_{i} C,$ defined by $\omega(\mathbf{x})= \widetilde{\omega}(\mathbf{x}) \mbox{ for any } \mathbf{x}\in \Gamma(G\ten_{\ell} H,u_{G}\ten_{\ell} u_{H}),$
is a homomorphism of  MV-algebras and $\omega(x\ten_{\ell} y)=\widetilde{\omega}(x\ten_{\ell} y)=(\widetilde{\omega}\circ\gamma_{G,H})(x,y)=\widetilde{\gamma}(x,y)=\gamma(x,y),$ for any $x\in \Gamma(G,u_{G})$ and $y\in\Gamma(H,u_{H})$.\\
In order to prove the uniqueness, let $\theta :\Gamma(G\ten_{\ell} H,u_{G}\ten_{\ell} u_{H})\rightarrow [0,\gamma(u_{G},u_{H})]\leq_{i} C$ be another  homomorphism of MV-algebras such that 
$\theta(x\ten_{\ell} y)=\gamma(x,y)$ for any $x\in \Gamma(G,u_{G})$ and $y\in\Gamma(H,u_{H})$. By \cite[Proposition 2.9]{LeuTens} there is a unique homomorphism of  $\ell$-groups $\widetilde{\theta}:G\ten H\rightarrow L$ such that 
$\widetilde{\theta}|_{\Gamma(G\ten_{\ell} H,u_{G}\ten_{\ell} u_{H})}=\theta$. If  $\tau$ is defined as $\widetilde{\theta}\circ\gamma_{G,H}$
then is straightforward that $\tau :G\times H\rightarrow L$ is a $\ell u$-bilinear function  and 
$\tau (x,y)=\theta (x\ten y)=\gamma(x,y)$. Since $\widetilde{\gamma }$ is the unique $\ell u$-bilinear function  that extends $\gamma$ we get $\tau=\widetilde{\gamma}$ and $\widetilde{\theta}\circ\gamma_{G,H}=\widetilde{\gamma}$.  It follows that  $\widetilde{\theta}=\widetilde{\omega}$ and $\theta=\omega$.
\end{proof}

In the following, we introduce a notation.
Let $A$, $B$ be MV-algebras and $(G_{A}, u_A)$, $(G_{B}, u_B)$ $\ell u$-groups such that 
$A\simeq \Gamma(G_{A}, u_A)$ and  $B\simeq \Gamma(G_{B}, u_B)$. Assume  $\{ \eta_{A} \}_{A\in \mathbf{MV}}$ is the natural isomorphism between the categories $\mathbf{MV}$ and $\mathbf{auG}$, i.e.   $\eta_A:A\to\Gamma(G_A,u_A)$ and 
 $\eta_B:B\to\Gamma(G_B,u_B)$ are isomorphisms of MV-algebra. Hence we define  
$\gamma_{A,B}:A\times B\to \Gamma (G_A\ten _{\ell} G_B, u_A \ten_{\ell} u_B)$ by 
$\gamma_{A,B}(x,y)= \eta_A(x)\ten_{\ell} \eta_B(y)$ for any $x\in A$ and $y\in B$.

\blem \label{teo:01un-l}
Let $A$, $B$  and $(G_{A}, u_A)$, $(G_{B}, u_B)$ $\ell u$-groups such that 
$A\simeq \Gamma(G_{A}, u_A)$ and  $B\simeq \Gamma(G_{B}, u_B)$.
For any MV-algebra $C$ and 
 any bimorphism $\gamma:A\times B\rightarrow C$ there is a unique homomorphism of MV-algebras
$\omega:\Gamma (G_A\ten _{\ell} G_B, u_A \ten_{\ell} u_B)\rightarrow [0,\gamma(1_{A},1_{B})]\leq_{i} C$ such that $\omega\circ\gamma_{A,B}=\gamma$.
\elem

\begin{proof}  
 In the following $A \ten_{\ell} B$ will denote the MV-algebra 
$ \Gamma (G_A\ten _{\ell} G_B, u_A \ten_{\ell} u_B)$.
Suppose that $C$ is an arbitrary MV-algebra and $\gamma:A\times B\rightarrow C$ is a bimorphism. We define $\gamma_{1}:\Gamma(G_{A}, u_A)\times \Gamma (G_{B},u_{B})\rightarrow \Gamma(G_{C},u_{c})$ by

$\gamma_{1}(x,y)=\eta_{C}(\gamma(\eta_{A}^{-1}(x),\eta_{B}^{-1}(y)))$.\\
Since $\eta_{C}$, $\eta_{A}^{-1}$ and $\eta_{B}^{-1}$ are
MV-algebra isomorphisms, $\gamma_{1}$ is also a bimorphism, and by Lemma \ref{lem:01ten}, 
there is a unique homomorphism of MV-algebras
$\omega_{1}:\Gamma(G_{A}\ten_{\ell} G_{B},u_{A}\ten_{\ell}u_{B})\rightarrow [0,\gamma_{1}(u_{A},u_{B})]\leq_{i} \Gamma(G_{C},u_{C})$ such that $\omega_{1}(x\ten_{\ell}y)=\gamma_{1}(x,y)$ 
for any $x\in \Gamma(G_{A},u_{A})$ and $y\in\Gamma(G_{B},u_{B})$.
Remark that the definition domain of $\omega_{1}$ is $A\ten_{\ell}B$. Hence, if we define $\omega(\mathbf{x})=\eta_{C}^{-1}(\omega_{1}(\mathbf{x}))\mbox{ for any  }\mathbf{x}\in A\ten_{\ell}B$, by \cite[Lemma 2.8]{LeuTens}  $\omega:A\ten_{\ell}B\rightarrow [0,\eta_{C}^{-1}(\gamma_{1}(u_{A},u_{B}))]\leq_{i}C$ is
 a homomorphism of MV-algebras. We have $\eta_{C}^{-1}(\gamma_{1}(u_{A},u_{B}))=\gamma(\eta_{A}^{-1}(u_{A}),\eta_{B}^{-1}(u_{B}))=\gamma(1_{A},1_{B}).$\\
For any $x\in A$ and $y\in B$ it follows that

$(\omega\circ\gamma_{A,B})(x,y)=\omega(\eta_{A}(x)\ten_{\ell}\eta_{B}(y))= \eta_{C}^{-1}(\omega_{1}(\eta_{A}(x)\ten_{\ell}\eta_{B}(y)))=$
 
$=\eta_{C}^{-1}(\gamma_{1}(\eta_{A}(x),\eta_{B}(y)))= \gamma(x,y)$.\\
Hence, $\omega\circ\gamma_{A,B}=\gamma$. In order to prove the uniqueness, suppose that
$\theta:A\ten_{\ell}B\rightarrow[0,\gamma(1_{A},1_{B})]\leq_{i} C$ is another homomorphism of MV-algebras 
 such that $\theta\circ\gamma_{A,B}=\gamma$. Using the isomorphism $\eta_{C}$ we define
 the following homomorphism of MV-algebras:
$\theta_{1}:A\ten_{\ell}B\rightarrow[0,\eta_{C}(\gamma(1_{A},1_{B}))]\leq \Gamma(G_{C},u_{C}),$ with $\theta_{1}(\mathbf{x})= \eta_{C}(\theta(\mathbf{x})),$ 
for any $\mathbf{x}\in A\ten_{\ell}B$. Remark that 

$\eta_{C}(\gamma(1_{A},1_{B}))=\eta_{C}(\gamma(\eta_{A}^{-1}(u_{A}),\eta_{B}^{-1}(u_{B})))=\gamma_{1}(u_{A},u_{B})$.\\
Moreover, for any $x\in \Gamma(G_{A},u_{A})$ and $y\in\Gamma(G_{B},u_{B})$ we get

$\theta_{1}(x\ten_{\ell}y)=\eta_{C}(\theta(x\ten_{\ell}y))=
\eta_{C}(\theta(\eta_{A}(\eta_{A}^{-1}(x))\ten_{\ell}\eta_{B}(\eta_{B}^{-1}(y))))=$

$=\eta_{C}(\theta(\gamma_{A,B}(\eta_{A}^{-1}(x),\eta_{B}^{-1}(y))))=\eta_{C}(\gamma(\eta_{A}^{-1}(x),\eta_{B}^{-1}(y)))=\gamma_{1}(x,y)$.\\
It follows that $\theta_{1}$ satisfies the properties that uniquely characterize $\omega_{1}$ by Lemma \ref{lem:01ten}, so $\theta_{1}=\omega_{1}$. In consequence, 

$\omega(\mathbf{x})=\eta_{C}^{-1}(\omega_{1}(\mathbf{x}))=\eta_{C}^{-1}(\theta_{1}(\mathbf{x}))=\theta(\mathbf{x})$\\
for any $\mathbf{x}\in A\ten_{\ell}B$, so $\theta=\omega$.
\end{proof}

We are now ready to  prove  the main result.

\bthm\label{pro:02ten-iso}
If $(G_A, u_A)$, $(G_B, u_B)$ are $\ell u$-groups and  $A$, $B$ are MV-algebras such that  $A\simeq \Gamma (G_A, u_A)$ and  $B\simeq \Gamma (G_B, u_B)$ then $A\ten_{MV}B\simeq  \Gamma (G_A \ten_{\ell}G_B, u_A \ten _{\ell}u_B)$.
\ethm 

\begin{proof}
As before,  $A \ten_{\ell} B$ denotes
$ \Gamma (G_A\ten _{\ell} G_B, u_A \ten_{\ell} u_B)$ and we prove that  $A\ten_{MV}B\simeq A\ten_{\ell}B$.
Let $\gamma_{A,B} : A\times B \rightarrow A \ten _{\ell} B$ the bimorphism defined as in Lemma \ref{teo:01un-l} and $\beta_{A,B}: A\times B\rightarrow A\ten_{MV}B$ the standard bimorphism of the MV-algebraic tensor product.
Using universal property of the tensor product, we get a homomorphism of  MV-algebras 

$\lambda : A\ten_{MV} B \rightarrow [0, \gamma_{A,B}(1_A, 1_B)]\le_i A\ten_{\ell} B$\\
such that $ \lambda_{A,B}\circ \beta_{A,B}= \gamma_{A,B}$. By Lemma \ref{teo:01un-l} there exists an homomorphism of MV-algebras 

$ \delta : A\ten_{\ell} B \rightarrow [0, \beta_{A,B}(1_A, 1_B)]\le_i A\ten_{MV}B$\\
such that $\delta \circ \gamma_{A,B}= \beta_{A,B}$.\\
Then we get:

$(\delta \circ \lambda )\circ \beta_{A,B}=\delta \circ (\lambda \circ \beta_{A,B})=\delta \circ \gamma_{A,B}=\beta_{A,B}$

$(\lambda \circ \delta )\circ \gamma_{A,B}=\lambda \circ (\delta \circ \gamma_{A,B})=\lambda  \circ \beta_{A,B}= \gamma_{A,B}$.\\
Therefore by the universal property of $\beta_{A,B}$ it follows $\delta \circ \lambda = \mathbf{I}_{A\ten_{MV}B}$, and by the universal property of $\gamma_{A,B}$ it follows $\lambda \circ \delta =\mathbf{I}_{A\ten_{\ell}B}$, that is the two tensor product are MV-algebraic isomorphic, i.e.  $\Gamma (G_A, u_A)\ten_{MV}\Gamma (G_B, u_B) \simeq \Gamma (G_A \ten_{\ell}G_B, u_A \ten _{\ell}u_B)$.
\end{proof}

\bfact 
We note that, in  \cite{LeuTens},  it  is proved that the functor $\Gamma$ commutes with another tensor product denoted $\ten_o$. 
The tensor product $\ten_o$ of MV-algebras defined in \cite{LeuTens} corresponds to the tensor product $\ten_o$ of $\ell$-groups defined in \cite{Mart} and it is defined using only bilinear functions, instead of  bimorphisms and $\ell u$-bilinear functions. In this case, if $A$ and $B$ are MV-algebras, the function $\iota_A:A\to A\ten_o B$ by $\iota_A(a)=a\ten_o 1$ is no longer a homomorphism of MV-algebras.
\efact

Recall that the functor $\Gamma$ maps archimedean $\ell u$-groups to semisimple MV-algebras.    Moreover, one can easily prove that $\Gamma$ also preserve the archimedean tensor product.

\bcor \label{isoTPss}
If $(G_A, u_A)$, $\Gamma (G_B, u_B)$ are archimedean  $\ell u$-groups and  $A$, $B$ are semisimple MV-algebras such that  $A\simeq \Gamma (G_A, u_A)$ and  $B\simeq \Gamma (G_B, u_B)$ then $A\ten B\simeq  \Gamma (G_A \ten_{a}G_B, u_A \ten _{a}u_B)$.
\ecor
\begin{proof} The proof is similar with the one of Theorem \ref{pro:02ten-iso}. The main idea is that $\Gamma (G_A \ten_{a}G_B, u_A \ten _{a}u_B)$ satisfy the same universal property that  uniquely defines $A\ten B$, up to isomorphism.
\end{proof}

\section{Scalar extension property for semisimple MV-algebras} \label{sec:SEP}

\noindent In this section we will investigate the scalar extension property for MV-algebras and $\ell u$-groups. First, let us  state the property in our context. 
\medskip

\noindent (SEP$_{MV}$) If $P$ is a unital  PMV-algebra and $A$ is an MV-algebra, then $P\ten_{MV}A$  has a canonical structure of MV-module over $P$.

\noindent (SEP$_{\ell}$) If $R$ is an $\ell$-ring and $G$ is an $\ell$-group, then $R\ten_{\ell}G$ has a canonical structure of $\ell$-module over $R$. 

\medskip

We summarize the results  so far.

\bfact\label{sep1} The scalar extension property  is one of the basic property arising from a tensor product, and while it is straightforward in the non-ordered case, with lattice ordered structures it presents some difficulties.\\
(1) For MV-algebras, the property SEP$_{MV}$ is stated in \cite[Theorem 4.11]{LeuTens}, but the proof presents a wrong argument. We note that the structure $P\ten_{MV}A$  can be endowed with a family of unary operations $\{\alpha\}_{\alpha\in P}$ such  that  the function $(\alpha, x)\mapsto \alpha x$ is linear in the second argument. The proof of the linearity in the first argument contains a mistake. \\
(2) For $\ell$-groups, the property SEP$_{\ell}$  is left as an exercise in \cite[Chapter 4.5]{Stein}. The case $R={\mathbb R}$, i.e. the $\ell$-modules are Riesz spaces (vector lattices),  is  considered in \cite[Proposition 2.1]{Mart}, but the details  are missing. \\
Our impossibility to correct the proof from \cite[Theorem 4.11]{LeuTens} and  to complete the  proof from \cite[Proposition 2.1]{Mart} is related to the fact that the sum of two homomorphisms of $\ell$-groups is not always an homomorphism of $\ell$-groups. We leave this as open problems. 
\efact

\bfact\label{sep2} Under the assumption that \cite[Proposition 2.1]{Mart} is correct, using Theorem \ref{pro:02ten-iso}, one can immediately prove (SEP$_{MV}$) for $P=[0,1]$. This means that for any MV-algebra $A$, the tensor product $[0,1]\ten_{MV}A$
is a Riesz MV-algebra. 
\efact

\bfact\label{sep3}
In \cite[Theorem 5]{BVanR}, the property SEP$_{\ell}$ is proved for the archimedean tensor product  in the context of Riesz spaces  ($R={\mathbb R}$). 
\efact

In Theorem \ref{teo:TPmod}  we prove SEP$_{MV}$   for the semisimple tensor product allowing $P$ to be an arbitrary unital and  semisimple PMV-algebra. Using  Theorem \ref{pro:02ten-iso} we get SEP$_{\ell}$  for archimedean structures allowing $R$ to be an 
arbitrary archimedean $\ell$-ring with strong unit. In this way we generalize the result from \cite{BVanR}.

\bfact \label{rem:TP01}
Let $P$ be a PMV-algebra and $S\subseteq P$ a sub PMV-algebra of $P$. It is easily seen that $P$ is a $S$-MV-module such that  the external operation coincides with the internal product on $P$.
\efact

\blem \label{lem:TP01}
Let $A\subseteq C(X)\subseteq C(X\times Y)$ be a PMV-algebra. Then the map 
\vsp{-0.2}
\begin{center}
$ \varphi :A\times C(X\times Y)\rightarrow C(X\times Y)$,

$ \varphi (a, f)(x,y)=a(x)f(x,y) $ for any $a\in A$ and $f\in C(X\times Y)$,
\end{center}
\vsp{-0.2}
defines a structure of $A$-MV-module.
\elem
\begin{proof}
If $a_1 + a_2$ is defined in $A$, then for any $x\in X$ the partial sum $a_1(x)+a_2(x)$ is defined and $(a_1+a_2)f=a_1f+a_2f$ if and only if the equality holds for any $x\in X, y\in Y$. Trivially  $(a_1(x)+a_2(x))f(x,y)=a_1(x)f(x,y)+a_2(x)f(x,y)$, since the functions are $[0,1]$-valued. In the same way we get all other conditions for an MV-module.
\end{proof}

\bthm \label{teo:TPmod}
Let $A$ be a unital and semisimple PMV-algebra, and $B$ be a semisimple MV-algebra. Then $A\ten B$ is an $A$-MV-module.
\ethm
\begin{proof}
By Theorem \ref{rem:TPMun}, $A\ten B= \left\langle \gamma (a,b) \mid a\in A \ b\in B \right\rangle _{MV} \subseteq C(X\times Y)$, with $A\subseteq C(X)$, $B\subseteq C(Y)$.\\
For any $\alpha \in A$, we define 

$\omega _{\alpha} : A\times B \rightarrow A\ten B$, \quad $ \omega_{\alpha}(a,b)(x,y)=\alpha (x) a(x)b(y).$\\
Since $a(x)b(y)=\gamma (a,b)(x,y)$ and $\alpha (x) \le 1$,  it follows that $\omega_{\alpha}(a,b) \in A\ten B$. We  prove that $\omega_\alpha$ is a bimorphism. 

Note that  $\omega_{\alpha}(a_1\wedge a_2, b)(x,y)=\alpha (x) [(a_1\wedge a_2)(x)]b(y)=\alpha (x) (a_1(x)\wedge a_2(x))b(y)$. Since we are in the unital case the product can be distributed on and we get $\wedge$ and $\vee$  $[\alpha (x)a_1(x) \wedge \alpha (x) a_2(x)]b(y)=[\alpha (x) a_1(x)b(y)] \wedge [\alpha (x)a_2(x)b(y)]$. In the same way we get the desired property for $\vee$ and on the second component;

If $a_1+a_2$ defined in $A$, then $\omega_{\alpha}(a_1+a_2, b)(x,y)=\alpha (x) [a_1(x)+a_2(x)]b(y)$. By the definition of PMV-algebra, we get the desired conclusion.\\
Therefore, applying the universal property,  there exists a homomorphism of MV-algebras

$ \Omega_{\alpha}:A\ten B \rightarrow [0, \omega_{\alpha}(1,1)]\le_i A\ten B .$\\
We remark that $\omega_{\alpha}(1,1)(x,y)=\alpha (x) \mathbf{1}(x) \mathbf{1}(y)= \alpha (x)$. Moreover, $\Omega_{\alpha}: A\ten B \rightarrow [0, \omega_{\alpha}(1,1)]\le_i A\ten B \subseteq C(X\times Y)$.\\
By Lemma \ref{lem:TP01}, $C(X\times Y)$ is $A$-MV-module with external operation $\varphi :A \times C(X\times Y)\rightarrow C(X\times Y)$. We fix $\alpha \in A$, and we define
\vsp{-0.2}
\begin{center}
$ \theta_{\alpha}: C(X\times Y)\rightarrow C(X\times Y), \quad \theta_{\alpha}(f)=\varphi (\alpha, f).$
\end{center}
\vsp{-0.2}
Since $C(X\times Y)$ is $A$-MV-module, $\theta_{\alpha}$ is linear and $\theta_{\alpha}(f) \le f$. It follows that $\theta_{\alpha}(f_1)\wedge f_2=0$ whenever $f_1\wedge f_2=0$. Note that such functions are called {\em \textit{f}-operators} in  \cite{LeuOp}.
 By \cite[Proposition 5.9]{LeuOp}, $\theta_{\alpha}$  is an MV-algebra morphism

$ \theta_{\alpha}: C(X\times Y)\rightarrow [0, \theta_{\alpha}(1)]\le_i C(X\times Y).$\\
Moreover, $\theta_{\alpha}(1)(x,y)=\varphi (\alpha, 1)(x,y)=\alpha (x)$.\\
Then we have

$ \Omega_{\alpha}: A\ten B \rightarrow [0, \omega_{\alpha}(1,1)]\le_i A\ten B\subseteq C(X\times Y)$ and 

$ \theta_{\alpha}|_{A\ten B}: A\ten B \rightarrow [0, \theta_{\alpha}(1)]\le_i C(X\times Y).$\\
Since $\omega_{\alpha}(1,1)=\theta_{\alpha}(1)$, by \cite[Lemma 2.3]{LeuTens},

$[0, \omega_{\alpha}(1,1)]_{A\ten B}\subseteq [0, \theta_{\alpha}(1)]\le_i C(X\times Y)$.\\
Therefore $\Omega_{\alpha}$ and $\theta_{\alpha}|_{A\ten B}$ are both maps from $A\ten B$ in $[0, \theta_{\alpha}(1)]\le_i C(X\times Y)$, and by universal property they are the same map if they coincide on generators, but this is trivially true, since for any $a\in A,\ b\in B,\ x\in X,\ y\in Y$

$\Omega_{\alpha}(\gamma (a,b))(x,y)=\alpha (x) a(x)b(y)= \varphi (\alpha, \gamma (a,b))(x,y)= \theta_{\alpha}(\gamma (a,b))(x,y)$.\\
As result, we have $A\ten B$ included in $C(X\times Y)$ as MV-algebra, and two families of linear functions:

$\{ \Omega_{\alpha}: A\ten B \rightarrow A\ten B \}_{\alpha \in A}$;\quad $\{ \theta_{\alpha}: C(X\times Y) \rightarrow C(X\times Y)\}_{\alpha \in A}$\\
such that $\theta_{\alpha}|_{A\ten B}=\Omega_{\alpha}$, and $(C(X\times Y), \{ \theta_{\alpha}\}_{\alpha \in A})$ is an $A$-MV-module.\\
Since for any $\alpha \in A$ we have $\theta_{\alpha}|_{A\ten B}(A\ten B) = \Omega_{\alpha}(A\ten B)\subseteq A\ten B$, $A\ten B$ is closed to the scalar product, and $(A\ten B, \{ \theta_{\alpha}|_{A\ten B} \}_{\alpha \in A})$ is an sub MV-module of $C(X\times Y).$ 
\end{proof}

By Theorem \ref{pro:02ten-iso} and Theorem \ref{teo:TPmod} we immediately get the following.

\bthm \label{teo:SEPgroups} If $R$ is a unital and archimedean $\ell u$-ring and $G$ is an archimedean $\ell u$-group, $R\otimes_{a} G$ is a $\ell u$-module over $R$.
\ethm

In the sequel, we apply use SEP$_{MV}$ in order to establish  connections between MV-algebras and Riesz MV-algebras.

\bprop \label{teo:TPriesz}
Let $A$ be a semisimple Riesz MV-algebra, and $B$ be a semisimple MV-algebra. Then $A\ten B$ is an Riesz MV-algebra. In particular, $[0,1]\ten B$ is a Riesz MV-algebra. 
\eprop
\begin{proof}
The proof is similar with the one of  Theorem \ref{teo:TPmod}. The starting bimorphism will be 

$\omega_{\alpha}: A \times B\to A\ten B$, $\omega_{\alpha}(a,b)=\alpha a(x) b(y)$ where $\alpha \in [0,1]$.
\end{proof}

By Theorem \ref{pro:02ten-iso} and Proposition \ref{teo:TPriesz} we immediately get the following.

\bprop \label{pro:SEPRiesz}
Let $V$ be an archimedean Riesz Space with strong unit and $G$ be an archimedean $\ell u$-group. Then $V\ten_a G$ is an archimedean Riesz Space with strong unit. In particular, $\mathbb{R}\ten_a G$ is a Riesz Space. 
\eprop

\bfact\label{eud} Let  $B$  be  a semisimple MV-algebra and  $\iota_B:B\to [0,1]\ten B$ the canonical embedding. Since $[0,1]\ten B$ is generated, as a Riesz MV-algebra, by $\iota_B(B)$ it follows by \cite[Corollary 4.2]{DDIL}, that $[0,1]\ten A$ is, up to isomorphism, the Riesz MV-algebra hull of $B$. One can see \cite{DDIL} for more details. 
\efact 

\bcor \label{cor:riesz}
Let $B$ be a  semisimple MV-algebra. For any  semisimple Riesz MV-algebra $V$ and for any homomorphism of MV-algebras $f:B\rightarrow \mathcal{U}_{\mathbb{R}}(V)$ there is a unique homomorphism of Riesz MV-algebras
$\widetilde{f}:[0,1]\ten B\rightarrow V$ such that $\widetilde{f}\circ \iota_{B}=f$.
\ecor
\begin{proof}
Define $\beta_f:[0,1]\times B\to V$ by $\beta(\alpha,x)=\alpha f(x)$ for any $\alpha\in [0,1]$ and $x\in B$ and use the universal property of the tensor product. We note that any MV-algebra homomorphism between Riesz MV-algebras preserves the scalar multiplication, so it is a morphism of Riesz MV-algebras \cite[Corollary 3.11]{LeuRMV}.
\end{proof}

Assume that $\mathbf{MV_{ss}}$ is the full subcategory of semisimple MV-algebras and $\mathbf{RMV_{ss}}$ is the full subcategory of semisimple Riesz MV-algebras. Hence we define a functor ${\mathcal T}_\ten:\mathbf{MV_{ss}}\to \mathbf{RMV_{ss}}$ by 

${\mathcal T}_\ten(B)=[0,1]\ten B$ for any semsimple MV-algebra $B$ and 

${\mathcal T}_\ten(f) =\widetilde{f}$ for any homomorphism of MV-algebras $f:A\to B$, where $\widetilde{f}:[0,1]\ten A\to [0,1]\ten B$ is the unique Riesz MV-algebra homomorphism such that $\widetilde{f}\circ \iota_A=\iota_B\circ f$ (which exists by  Corollary \ref{cor:riesz}).

\bcor \label{teo:adjMVRiesz}
Under the above hypothesis,  ($\mathcal{T}_\ten, \mathcal{U}_{\mathbb R})$  is an adjoint pair. 
\ecor
\begin{proof} It is straightforward.   Using a different construction of the Riesz hull, this result is proved in \cite{DDIL}.
\end{proof}

\section{The categorical adjunction between semisimple PMV-algebras and semisimple \textit{f}MV-algebras} \label{sec:PMV}

In the sequel, using the scalar extension property, we define an adjunction between the category of semisimple and unital  PMV-algebras and the category of unital and semisimple \textit{f}MV-algebras. 

As a preliminary step, we   prove  the following theorem.

\bprop\label{teo:TPprod}
Let $A$, $B$ be unital and semisimple PMV-algebras. Then $A\ten B$ is a unital and semisimple PMV-algebra.
\eprop
\begin{proof}
By Theorem \ref{rem:TPMun}, $A\ten B= \left\langle \gamma (a,b) \mid a\in A \ b\in B \right\rangle _{MV} \subseteq C(X\times Y)$ with $A\subseteq C(X)$, $B\subseteq C(Y)$, as MV-algebra.\\
For any $c \in A\ten B$, we define 
$\omega _{c}:A\times B \rightarrow A\ten B$,  $ \omega_{c}(a,b)(x,y)=c (x,y) a(x)b(y).$\\
Since $a(x)b(y)=\gamma (a,b)(x,y)$ and $c(x,y) \le 1$, $\omega_{\alpha}(a,b) \in A\ten B$ and it is a bimorphism, likewise in the proof of Theorem \ref{teo:TPmod}. Therefore, applying the universal property, there exist a map

$ \Omega_{c}:A\ten B \rightarrow [0, \omega_{c}(1,1)]\le_i A\ten B .$\\
Again, $\omega_{c}(1,1)(x,y)=c (x,y) \mathbf{1}(x) \mathbf{1}(y)= c (x,y)$ and $\Omega_{c}: A\ten B \rightarrow [0, \omega_{c}(1,1)]\le_i A\ten B \subseteq C(X\times Y)$.\\
It is straightforward that $C(X\times Y)$ is PMV-algebra with internal product $*:C(X\times Y) \times C(X\times Y)\rightarrow C(X\times Y)$ defined component-wise. We fix $c \in A\ten B$, and we define

$ \theta_{c}: C(X\times X)\rightarrow C(X\times Y), \quad \theta_{c}(f)=c*f.$\\
It is easy to prove, since $C(X\times Y)$ is unital PMV-algebra, that $\theta_{c}$ is linear and $\theta_{c}(f) \le f$. And again like in the proof of Theorem \ref{teo:TPmod}, $ \theta_{c}: C(X\times Y)\rightarrow [0, \theta_{c}(1)]\le_i C(X\times Y)$ is an homomorphism of MV-algebras. Moreover, $\theta_{c}(1)(x,y)=(c*1)(x,y)=c(x,y)$.\\
Then we have

$ \Omega_{c}: A\ten B \rightarrow [0, \omega_{c}(1,1)]\le_i A\ten B\subseteq C(X\times Y)$ and 

$ \theta_{c}|_{A\ten B}: A\ten B \rightarrow [0, \theta_{c}(1)]\le_i C(X\times Y).$\\
Since $\omega_{c}(1,1)=\theta_{c}(1)$, we get $[0, \omega_{c}(1,1)]_{A\ten B}\subseteq [0, \theta_{c}(1)]\le_i C(X\times Y)$, and the conclusion follows like in Theorem \ref{teo:TPmod}.\\
Finally, it can be easily seen that the unit in $A\ten B$ is the element $1_A \ten 1_B$. Then $A\ten B$ is unital, therefore it is semisimple as PMV-algebra.
\end{proof}

\bthm \label{pro:TPfMV}
Let $R$ be a \textit{f}MV-algebra and $P$ be a unital and semisimple PMV-algebra. Then $R\ten P$ is a unital and semisimple \textit{f}MV-algebra.
\ethm
\begin{proof}
By construction $R\ten P$ is a semisimple MV-algebra, by Corollary \ref{teo:TPriesz}, $R\ten P$ is a Riesz MV-algebra and by Proposition \ref{teo:TPprod} it is a unital and semisimple PMV-algebra. Moreover, by the construction of the product and the scalar operation as the usual product and scalar operation between functions given in Proposition \ref{teo:TPprod} and Theorem \ref{teo:TPmod}, the associativity law between products is satisfied since it holds in any $C(X)$. It follows that $R\ten P$ is a unital and semisimple \textit{f}MV-algebra.
\end{proof}

\bprop \label{prop:adj}
Let $A$ be a unital and semisimple PMV-algebra. For any  unital and semisimple \textit{f}MV-algebra $M$ and for any homomorphism of PMV-algebras $f:A\rightarrow \mathcal{U}_{\mathbb{R}}(M)$ there is a unique homomorphism of \textit{f}MV-algebras
$\widetilde{f}:[0,1]\ten A\rightarrow M$ such that $\widetilde{f}\circ \iota_{A}=f$, where $\iota_{A} : A \rightarrow [0,1]\ten A$ is the embedding in the tensor product.
\eprop
\begin{proof}
By Corollary \ref{cor:riesz}, there exists a homomorphism of Riesz MV-algebras $\widetilde{f}:[0,1]\ten A \rightarrow \mathcal{U}_{\mathbb{(\cdot)}}(M)$ such that $\widetilde{f} \circ \iota_A = f$. Since $[0,1]\ten A$ is a \textit{f}MV-algebra by Theorem \ref{pro:TPfMV}, $\widetilde{f}$ is a homomorphism of Riesz MV-algebras between unital and semisimple \textit{f}MV-algebras. By \cite[Proposition 3.2]{LL1} $\widetilde{f}$ is a homomorphism of \textit{f}MV-algebras.
\end{proof}

\bprop \label{pro:TPhom}
Let $P_1$ and $P_2$ be semisimple and unital PMV-algebras, and $h:P_1 \to P_2$ a homomorphism of PMV-algebras. Then there exists a unique $h^{\sharp}:[0,1]\ten P_1 \to [0,1]\ten P_2$ homomorphism of \textit{f}MV-algebras such that $h^{\sharp} \circ \iota_1 = \iota_2 \circ h$, where $\iota_i:P_i \to [0,1]\ten P_i$ for $i=1,2$ are the natural embeddings.
\eprop
\begin{proof}
It is straightforward by Proposition \ref{prop:adj}, with $f=\iota_2\circ h$.
\end{proof}

\noindent Let $\mathbf{uPMV_{ss}}$ be the full subcategory of semisimple and unital PMV-algebras with homomorphism of PMV-algebras  and let  $\mathbf{ufMV_{ss}}$ be the full  subcategory of semisimple and unital \textit{f}MV-algebras with homomorphism of  \textit{f}MV-algebras.\\
We define a functor $\mathcal{F}_{\ten}: \mathbf{uPMV_{ss}} \rightarrow \mathbf{ufMV_{ss}}$ as follows
\begin{itemize}
\item[(i)] for any $P\in \mathbf{uPMV_{ss}}$, $\mathcal{F}_{\ten}(P)$ is $[0,1]\ten P$. By Theorem \ref{pro:TPfMV} it is a unital, commutative and semisimple \textit{f}MV-algebra.
\item[(ii)] for any homomorphism of  PMV-algebras $h:P_1\rightarrow P_2$, $\mathcal{T}(h)$  is the homomorphism of  \textit{f}MV-algebras $h^{\sharp}$ defined in Proposition \ref{pro:TPhom}.
\end{itemize}
\noindent From $\mathbf{ufMV_{ss}}$ to $\mathbf{uPMV_{ss}}$ we have the usual forgetful functor $\mathcal{U}_{\mathbb{R}}$.

\blem
$\mathcal{F}_{\ten}$ is a functor.
\elem
\begin{proof}
Let $P_1,P_2,P_3$ be PMV-algebras and $h:P_1\rightarrow P_2$ and $g:P_2 \rightarrow P_3$ PMV-algebras homomorphisms. Let $\iota_1 , \iota_2 $ and $\iota_3$ be the embeddings of $P_1, P_2, P_3$ in $[0,1]\ten P_1$, $[0,1]\ten P_2$ and $[0,1]\ten P_3$ respectively. Then 
\vsp{-0.2}
\begin{center}
$(g^{\sharp}\circ h^{\sharp}) \circ \iota_1 = g^{\sharp} \circ (h^{\sharp}\circ \iota_1)= g^{\sharp} \circ (\iota_2 \circ h)= (g^{\sharp} \circ \iota_2) \circ h = \iota_3 \circ (g\circ h)$
\end{center}
\vsp{-0.2}
and by Proposition \ref{pro:TPhom} we get $g^{\sharp}\circ h^{\sharp}= (g\circ h)^{\sharp}$.
\end{proof}
\blem
The maps $\{ \iota_A \}_{A\in \mathbf{uPMV_{ss}}}$ are a natural transformation between the identity functor on $\mathbf{uPMV_{ss}}$ and $\mathbf{ufMV_{ss}}$.
\elem
\begin{proof}
Let $P_1,P_2 \in \mathbf{uPMV_{ss}}$ and let $h:P_1\rightarrow P_2$ be a homomorphism of PMV-algebras. We need to prove that $\mathcal{U}_{\mathbb R}\mathcal{F}_{\ten}(h) \circ \iota_1 = \iota_2 \circ h$. Since $\mathcal{U}_{\mathbb R}\mathcal{F}_{\ten}(h)=h^{\sharp}$ the result follows from Proposition \ref{pro:TPhom}.
\end{proof}

\bthm \label{teo:adjPMVfMV}
The functors $\mathcal{F}_{\ten}$ and $\mathcal{U}_{\mathbb R}$ are adjoint functors.
\ethm
\begin{proof}
In order to prove that $\mathcal{F}_{\ten}$ is left adjoint functor of $\mathcal{U}_{\mathbb R}$, we need to prove that for any semisimple and unital \textit{f}MV-algebra $A$ and any homomorphism of PMV-algebras $f: P\rightarrow \mathcal{U}_{\mathbb R}(A)$, with $P\in \mathbf{uPMV_{ss}}$, there exists a homomorphism of \textit{f}MV-algebras $f^{\sharp}: \mathcal{T}(P)\rightarrow A$ such that $\mathcal{U}_{\mathbb R}(f^{\sharp})\circ \iota_P = f$. This follows from Proposition \ref{prop:adj}.
\end{proof}

Results in the previous sections can transferred to $\ell u$-groups and all related structures. We remark that Proposition \ref{teo:TPprod}, Theorem \ref{pro:TPfMV} and Corollary \ref{pro:02ten-iso} entail the following.

\bprop \label{pro:SEPfMV}

(i) If $R$ and $S$ are unital and archimedean $\ell u$-rings, $R\otimes_a S$ is a unital and archimedean $\ell u$-ring.

(ii) If $V$ is a unital and archimedean \textit{fu}-algebra and $R$ is a unital and archimedean $\ell u$-ring, $V\otimes_a R$ is a unital and archimedean \textit{fu}-algebra.
\eprop
\begin{proof}
It is straightforward by Theorem \ref{pro:TPfMV}, Corollary \ref{pro:02ten-iso} and the categorical equivalence. 
\end{proof}

\section{Conclusions}\label{sec:final}

By categorical equivalence, the adjunctions $(\mathcal{T}_\ten , \mathcal{U}_{\mathbb{R}})$ and $(\mathcal{F}_\ten , \mathcal{U}_{\mathbb{R}})$ naturally transfer to lattice-ordered structure. We denote by $\mathbf{auG_a}$ the category of archimedean $\ell u$-groups; $\mathbf{uR_a}$ the category of archimedean and unital $\ell u$-rings; $\mathbf{uRS_a}$ the category of archimedean Riesz Spaces with strong unit; $\mathbf{fuAlg_a}$ the category of archimedean and unital \textit{fu}-algebras.

Applying the inverse of  $\Gamma$ and $\Gamma_{\mathbb R}$, $(\mathcal{T}_\ten , \mathcal{U}_{\mathbb{R}})$ extends to $(\mathcal{T}_{\ten a} , \mathcal{U}_{\ell \mathbb{R}})$. This is an adjunction between $\mathbf{auG_a}$ and $\mathbf{uRS_a}$.

Applying the converses of the functors $\Gamma_{(\cdot )}$ and $\Gamma_{f}$ , $(\mathcal{F}_\ten , \mathcal{U}_{\mathbb{R}})$ extends to $(\mathcal{F}_{\ten a} , \mathcal{U}_{\ell \mathbb{R}})$. This is an adjunction between $\mathbf{uR_a}$ and $\mathbf{fuAlg_a}$.

In the following we summarize or results.

\bthm
 The following diagrams are commutative:
\begin{center}
\begin{tikzpicture}
  \node (A) {$\mathbf{uR_a}$};
  \node (B) [below of=A] {$\mathbf{PMV_{ss}}$};
  \node (C) [right of=A] {$\mathbf{auG_a}$};
  \node (D) [below of=C] {$\mathbf{MV_{ss}}$};
  \node (E) [right of=C] {$\mathbf{uRS_a}$};
  \node (F) [below of=E] {$\mathbf{RMV_{ss}}$};
\node(H)[left of=A] {$\mathbf{fuAlg_a}$};
\node(G)[below of=H] {$\mathbf{fMV_{ss}}$};

 \draw[->] (A) to node [swap] {$\Gamma_{(\cdot)}$} (B);
\draw[->] (H) to node [swap] {$\Gamma_{f}$} (G);
   \draw[->] (C) to node [swap] {$\Gamma$} (D);
\draw[->] (E) to node {$\Gamma_{\mathbb{R}}$} (F);

\draw[->,bend right] (A) to node [swap]{$\mathcal{F}_{\ten_a}$} (H);  
\draw[->] (H) to node [swap]{$\mathcal{U}_{( \ell \mathbb{R})}$} (A);  
\draw[->] (G) to node {$\mathcal{U}_{\mathbb{R}}$} (B);
\draw[->,bend left] (B) to node {$\mathcal{F}_{\ten}$} (G);  

\draw[->,bend left] (C) to node {$\mathcal{T}_{\ten_a}$} (E);  
\draw[->] (E) to node {$\mathcal{U}_{( \ell \mathbb{R})}$} (C);  
\draw[->] (F) to node [swap]{$\mathcal{U}_{\mathbb{R}}$} (D);
\draw[->,bend right] (D) to node [swap]{$\mathcal{T}_{\ten}$} (F);

\end{tikzpicture}
\end{center}
\ethm
\begin{proof} It is a direct consequence of Theorem \ref{pro:02ten-iso}. \end{proof}

Recall that  MV-algebras  were defined as the algebraic structures corresponding to \L ukasiewicz $\infty$-valued logic. Even if their algebraic theory  is relevant in itself as proved by \cite{CDM, MunBook}, it has been developed in strong connection with the associated logical system. The same holds for PMV-algebras \cite{MonPMV}, their theory has its origins in the problem of enriching
 \L ukasiewicz logic with a  binary conector whose interpretation in $[0,1]$ is the natural product. Logical systems were also defined for Riesz MV-algebras \cite{LeuRMV} and \textit{f}MV-algebras \cite{LL1}. Note that for PMV-algebras and \textit{f}MV-algebras the logical systems are developed only for particular suitable subclasses.   One important link between logic and algebra  in all these cases is the  Lindenbaum-Tarski  algebra, which is the free algebra generated by the set of propositional variables. Due to the fact  that the  free algebras are semisimple structures,  the scalar extension property allows us to connect the free structures as in Proposition \ref{last}.

For a nonempty set $X$, let $Free_{MV}(X)$  and  $Free_{RMV}(X)$ be the free MV-algebra and, respectively, the free Riesz MV-algebras generated by $X$. Let $Free_{PMV}(X)$  be the free PMV-algebra in $HSP([0,1]_{PMV})$, the variety of PMV-algebras generated by $[0,1]$. Similarly, let $Free_{\textit{f}MV}(X)$  be the free \textit{f}MV-algebra in $HSP([0,1]_{fMV})$, the variety of \textit{f}MV-algebras generated by $[0,1]$. See more details in \cite{CDM, MonPMV, LeuRMV, LL1}.

\bprop\label{last} For any nonempty set $X$, the following hold:\\
(i) $Free_{RMV}(X)\simeq [0,1]\ten Free_{MV}(X)$,\\
(ii) $Free_{fMV}(X)\simeq [0,1]\ten Free_{PMV}(X)$.
\eprop
\begin{proof}
(i) follows by Remark \ref{eud} and \cite[Proposition 4.1]{DDIL}; it can also be proved directly, similarly with (ii).\\
(ii) Assume $V$ is a \textit{f}MV-algebra and  $f:X\to V$ is  a function. Hence there is a unique homomorphism of PMV-algebras 
$f^{\#}: Free_{PMV}(X)\to {\cal U}_{\mathbb R}(V)$  which extends $f$.  By Proposition \ref{prop:adj}, there exists a homomorphism of \textit{f}MV-algebras $\widetilde{f}:[0,1]\ten Free_{PMV}(X)\to V$ such that $\widetilde{f}\circ\iota_{Free_{PMV(X)}}=f^{\#}$, so $\widetilde{f}(1\ten x)=f(x)$ for any $x\in X$. The uniqueness of $\widetilde{f}$ is a consequence of the uniqueness of $f^\#$. Since $\iota_{Free_{PMV(X)}}$ is an embedding we have 
$X\simeq \{1\ten x\mid x\in X\}$ so $[0,1]\ten Free_{PMV}(X)$ satisfies the universal property that uniquely defines $Free_{fMV}(X)$.
\end{proof}


\begin{thebibliography}{}
\bibitem{LeuGS}
Ball R.N., Georgescu G., Leu\c stean I., \textit{Cauchy completions of MV-algebras}, Algebra Universalis 47(4) (2002) 367-407.

\bibitem{BKW}
Bigard A., Keimel K.,   Wolfenstein S., \textit{ Groupes et anneaux r\'{e}ticul\'{e}s}, Lectures Notes in Mathematics Vol. 608 (1977), Springer-Verlag. 

\bibitem{Birk}
Birkhoff G., \textit{Lattice Theory}, AMS Coll. Publ. 25 3rd Ed. 1973.

\bibitem{BP}
Birkhoff G., Pierce R.S., \textit{Lattice-ordered rings}, An. Acad. Brasil. Cienc. 28 (1956) 41-69.

\bibitem{BVanR}
Buskes G.J.H.M., Van Rooij A.C.M., \textit{The Archimedean $\ell$-group tensor product}, Order 10 (1993) 93-102.

\bibitem{Cha1}
Chang C.C., \textit{Algebraic analysis of many valued logics}, Trans. Amer. Math. Soc. 88 (1958) 467-490.

\bibitem{Cha2}
Chang C.C.,  \textit{A new proof of the completeness of the \L ukasiewicz axioms}, Transactions of the American Mathematical Society 93 (1959),  74-80.

\bibitem{CDM}
Cignoli R., D'Ottaviano I.M.L., Mundici D., \textit{Algebraic foundation of many valued Reasoning}, Kluver Academc Publ Dordrecht 2000.


\bibitem{DiND}
Di Nola A., Dvurecenskij A., \textit{Product MV-algebras}, Multiple-Valued Logics 6 (2001) 193-215.

\bibitem{LeuMod}
Di Nola A., Flondor P., Leustean I., \textit{MV-modules}, Journal of Algebra 267(1) (2003) 21-40.

\bibitem{LeuRMV}
Di Nola A., Leustean I., \textit{\L ukasiewicz logic and Riesz Spaces}, Soft Comp. , to appear.

\bibitem{DDIL}
Diaconescu D., Leu\c stean I., \textit{The Riesz hull of a semisimple MV-algebra}, Mathematica Slovaca (special issue in honor of Antonio Di Nola), in print.

\bibitem{MoFl}  
Flaminio T., Montagna F., private communication.

\bibitem{LeuTens}
Flondor P., Leu\c stean I., \textit{Tensor Product of MV-algebras}, Soft Computing 7 (2003) 446-457.

\bibitem{LeuOp}
Flondor P., Leu\c stean I., \textit{MV-algebras with operator (the commutative and non-commutative case)}, Discrete Mathematics 274 (2004) 41-76.

\bibitem{Fre}
Fremlin D.H., \textit{Tensor Product of Archimedean Vector Lattices}, American Journal of Math. 94(3) (1972) 777-798.

\bibitem{LL1}
Lapenta S., Leu\c stean I., \textit{Towards Pierce-Birkhoff conjecture via MV-algebras}, submitted.

\bibitem{Mart}
Martinez J., \textit{Tensor Product of partially ordered groups}, Pac. J. Math. 41 (1972) 771-789.
 
\bibitem{MonPMV}
Montagna F., \textit{An algebraic approach to Propositional Fuzzy Logic}, Journal of Logic, Language and Information 9 (2000) 91-124.

\bibitem{Mun1}
Mundici D., \textit{Interpretation of ACF*-algebras in \L ukasiewicz sentential calculus}, J. Funct. Anal. 65 (1986) 15-63.

\bibitem{Mun}
Mundici D., \textit{Tensor products and the Loomis-Sikorski theorem for MV-algebras}, Advanced in Applied Mathematics 22 (1999) 227-248.

\bibitem{MunBook}
Mundici D., \textit{Advances in \L ukasiewicz calculus and MV-algebras}, Trends in Logic 35 Springer 2011.

\bibitem{Stein}
Steinberg S.A.,  \textit{Lattice-Ordered Rings and Modules}, Springer, 2010. 

\bibitem{RSZan}
Zaneen A.C., \textit{Riesz Space II}, North Holland, Amsterdam 1983.


\end{thebibliography}
\end{document}